\magnification=1200  %=`magstep 0.5
\input mssymb
\outer\def\beginsection#1\par{\filbreak\bigskip\leftline{\bf#1}\nobreak
\smallskip\vskip-\parskip\noindent}
\def\spec#1#2{\hbox{\dimen0=#1truein\vrule
\vbox{\tabskip= \dimen0\baselineskip= 20pt
\halign{\kern -.9\dimen0\rlap{$\scriptstyle##$}\hfil&&\rlap
{$\scriptstyle##$}\hfil\cr#2\crcr}\smallskip\hrule}}}
\def\aut{{\rm Aut}}
\def\tilp{{\widetilde P}}
 % the point at start of n-z spec entries
 % the indicator for zero spec entries
\def\mapright#1{\smash{\mathop{\longrightarrow}\limits^{#1}}}
\def\mapdown#1{\Big\downarrow\rlap{$\vcenter{\hbox{$\scriptstyle#1$}}$}}
\def\proof{\medbreak \noindent{\it Proof. }}
\def\sone {\ifmmode {{\rm S}^1} \else$\sone$\fi}
\def\beginsubsection #1. {\medbreak \noindent{\bf #1.}\enspace }

\def\res {{\rm Res}}
\def\cor {{\rm Cor}}
\def\hom {{\rm Hom}}
\def\pra{\par}
%from this point on are definitions added to ch2

\def\coh #1#2#3{H^{#1}(#2;#3)}
\def\cohz #1#2{\coh {#1} {#2} {\Bbb Z}}
\def\cohp #1#2{H^{#1}({#2};\Bbb F_p)}

\def\co #1#2{H^{#1}(#2)}
\def\chain{\dot{\hbox{\kern0.3em}}}
\def\pone{{\rm P^1}}%temporary definition
\def\qed{\hfill${\scriptstyle\blacksquare}$}

\def\class{{\rm B}}
\def\bg{{\class G}}
\def\bq{{\class Q}}
\def\bpsymb #1{{\class P(#1)}}
\def\bpn{\bpsymb n}
\def\bpt{\bpsymb 3}
\def\bpf{\bpsymb 4}
\def\tilg {{\widetilde G}}

\def\book#1/#2/#3/#4/#5/#6/ {\item{[#1]} #2{\sevenrm #3}. {\it #4}.\
#5, (#6).\par\smallskip}
\def\paper#1/#2/#3/#4/#5/#6/#7/#8/ {\item{[#1]} #2{\sevenrm #3}. #4.
{\it #5} {\bf #6} (#7), #8.\par\smallskip}
\def\prepaper#1/#2/#3/#4/#5/ {\item{[#1]} #2{\sevenrm #3}. #4.
#5\par\smallskip}
%definitions for references beyond here.
\def\alzub{1}
\def\araki{2}
\def\atiyah{3}
\def\brown{4}
\def\hueb{5}
\def\huec{6}
\def\leary{7}
\def\lew{8}
\def\moselle{9}
\def\thompap{10}
\def\Thomas{11}
\def\vasquez{12}

\baselineskip= 20 pt

\centerline{The integral cohomology rings of some $p$-groups}
\centerline{I.~J.~Leary}
\centerline{Trinity College, Cambridge} \par\bigskip
\beginsection 1. Introduction.\par
We determine the integral cohomology rings of an infinite family of $p$-groups,
for odd primes $p$,
with cyclic derived subgroups. Our method involves embedding the
groups in a compact Lie group of dimension one, and was suggested independently
by P.~H.~Kropholler and J.~Huebschmann. This construction has also been used by
the author to calculate the mod-$p$ cohomology of the same groups and by
B.~Moselle to obtain partial results concerning the mod-$p$ cohomology of
the extra special $p$-groups [\leary], [\moselle].
\par

\beginsection 2 The method and the groups.\par
Given a finite group $G$ and a central cyclic subgroup $C$, we fix an
embedding of $C$ into \sone, and define a Lie group $\tilg$ as the product of
\sone\ and $G$ amalgamating $C$, that is
$$\tilg={{\sone\times G}/\{(c^{-1},c):c\in C\}}$$
Then we have a commutative diagram:
$$\matrix{C & \mapright{} & G & \mapright{} & Q \cr \mapdown{} && \mapdown{} &
& \mapdown{} \cr
\sone & \mapright{} & \widetilde G & \mapright{} & Q.} $$
If $M$ is a $G$-module on which $C$ acts trivially, we may consider $M$ as a
$\widetilde G$-module by letting \sone\ act trivially, and the
Lyndon-Hochschild-Serre spectral sequence for the second extension
is often simpler than that for the first. To find $H^*({\rm B}G;M)$, given
$H^*({\rm B}\widetilde G;M)$, we use the Serre spectral sequence
of the fibration
$$\sone /C \cong {\widetilde G / G} \mapright{} {\rm B}G \mapright{}
{\rm B}\widetilde G.$$
This spectral sequence has $E_2^{i,j}=0$ for $j>1$, so the only possible
non-zero differential is $d_2$.
The above was first suggested to the author by P.~Kropholler.
A similar idea occurs in J.~Huebschmann's papers
[\hueb] and [\huec]. In the case where $M$ is a commutative
ring on which $G$ acts trivially, it appears that we may obtain another
filtration of $H^*({\rm B}G;M)$ by examining
the Eilenberg-Moore spectral sequence for the pullback square:
$$\matrix{{\rm B}G&\mapright{}&{\rm B}\widetilde G\cr
\mapdown{}&&\mapdown{}\cr
\{*\}&\mapright{}&{\rm B}{\widetilde G/G},}$$
but it can be shown that the two filtrations are identical. These spectral
sequences are just alternative ways to view the Gysin sequence for the
\sone-bundle ${\rm B}G$ over ${\rm B}\tilg$.
\par
It is possible for non-equivalent extensions of $C$ by $Q$ to yield equivalent
extensions of \sone\ by~$Q$. In fact this happens if and only if their
extension classes in $\coh 2 Q C$ map to the same element of $\coh 3 Q {\Bbb Z}$
under the Bockstein associated with the coefficient sequence $\Bbb
Z\rightarrowtail \Bbb Z \twoheadrightarrow C$.
\par
The groups we shall consider are central extensions of $C_{p^{n-2}}$ by~
$C_p\oplus C_p$ where $p$ is an odd prime, and $n\geq 3$. They may be presented
as $$P(n)=\langle A,B,C\mid A^p=B^p=C^{p^{n-2}}=[A,C]=[B,C]=1,
[A,B]=C^{p^{n-3}}\rangle.$$
We shall let $\tilp$ be the corresponding central extension of \sone\ by~
$C_p\oplus C_p$, that is the group obtained from $\sone\times P(n)$ by
amalgamating the subgroup of $P(n)$ generated by $\langle C\rangle$ and the
$C_{p^{n-2}}$ subgroup of \sone. There are four central extensions of
$C_{p^{n-2}}$ by $C_p\oplus C_p$; two abelian ones, $P(n)$, and a metacyclic
group $M(n)$ containing a cyclic subgroup of index $p$. This may be checked by
verifying that the action of $\aut (C_p\oplus C_p)$ on $\coh 2
{{\rm B}(C_p\oplus C_p)}
{C_{p^{n-2}}}$ has only four orbits, and then explicitly constructing four
non-isomorphic groups. There are however only two central extensions of \sone\
by $C_p\oplus C_p$; the direct product, which is abelian, and $\tilp$ which is
not. This follows from the fact that $\aut (C_p\oplus C_p)$ acts transitively
on the non-zero elements of $\coh 3 {{\rm B}(C_p\oplus C_p)} {\Bbb Z}$, which
may be identified with $\coh 2 {{\rm B}(C_p\oplus C_p)} \sone$ via the
Bockstein associated with the coefficient sequence $\Bbb Z\rightarrowtail \Bbb R
\twoheadrightarrow \sone$. Hence we see that ${\rm B}M(n)$ is also an
\sone-bundle over ${\rm B}\tilp$,
and in fact $\coh * {{\rm B}M(n)} {\Bbb Z}$ could
easily be determined from the results of this paper. This cohomology ring has
already been calculated using other methods [\Thomas].
\beginsection 3. Calculations.\par
We now begin our calculation of $\co * {{\rm B}\tilp}$ by examining the
spectral sequence with integer coefficients for $\tilp$ considered as an
extension of \sone\ by $C_p\oplus C_p$.
The $E_2$ page is readily seen to be generated by elements $\alpha, \beta\in
E_2^{2,0}$, $\gamma\in E_2^{3,0}$ and $\tau \in E_2^{0,2}$ subject only to the
relations $p\alpha=p\beta=0$, $p\gamma=0$ and $\gamma^2=0$. Note that $\tau $
has infinite order.
Since $E_2^{i,j}$ is trivial for $j$ odd, we see that all the even
differentials must vanish. The behaviour of the differentials is summarised in
the following lemma.
\proclaim Lemma 1. In the above spectral sequence there are exactly two
non-zero differentials, $d_3$ and $d_{2p-1}$. $d_3(\tau)$ is a non-zero
multiple of $\gamma$, and $E_4$ is generated by the classes of the elements
$\alpha, \beta, p\tau,\ldots, p\tau^{p-1}, \tau^p$ and $\tau^{p-1}\gamma$. All of
these generators are universal cycles except for $\tau^{p-1}\gamma$, which is
mapped by $d_{2p-1}$ to a non-zero multiple of $\alpha^p\beta-\beta^p\alpha$.
The $E_\infty$ page is generated by the elements $\alpha, \beta, p\tau,\ldots,
p\tau^{p-1}, \tau^p$ subject only to the relations they satisfy as elements of
$E_2$, and  the relation $\alpha^p\beta =\beta^p\alpha$.
\par
\proof The derived subgroup of $\tilp$ consists of the subgroup of
its central $\sone $ of
order $p$, so there can be no homomorphism from $\tilp$ to $\sone$
that restricts to an isomorphism from the centre to $\sone$. It follows by
considering the natural isomorphism $\cohz 2 \bg \cong \hom(\bg,\sone)$ that
the element $\tau$ cannot survive to $E_\infty$, so we must have $d_3(\tau)$ a
non-zero multiple of $\gamma$. This determines $d_3$ completely.
It may be checked that $E_4$ is isomorphic to
the subring of $E_2$ generated by
$\alpha, \beta, p\tau,\ldots, p\tau^{p-1}, \tau^p$ and $\tau^{p-1}\gamma$. All
these elements must be universal cycles, with the possible exception of
$\tau^{p-1}\gamma$, because the groups in which their images under $d_n$ lie
are already trivial. The only remaining potentially non-zero differential is
$d_{2p-1}(\tau^{p-1}\gamma)$. To complete this proof it suffices to show that
in the $E_\infty $ page the relation $\alpha^p\beta=\beta^p\alpha$ must
hold.\par
Let $Q$ be the quotient of $\tilp$ by its \sone\ subgroup, and take generators
$\alpha',\beta'$ for $\cohz 2 \bq $ and $\gamma'$ for $\cohz 3 \bq $. The
statement that $\gamma $ does not survive to $E_\infty $ in the spectral
sequence is equivalent to the statement that $\gamma'$ is mapped to zero by the
inflation map from $Q$ to $\tilp$. Now we calculate $\phi(\gamma')$, where
$\phi$ is the integral cohomology operation $\delta_p\pone\pi_*$, where $\pi_*$
is the map induced by the change of coefficients from $\Bbb Z$ to $\Bbb F_p$,
$\pone$ is a reduced power, and $\delta_p$ is the Bockstein for the sequence
$\Bbb Z \rightarrowtail \Bbb Z \twoheadrightarrow \Bbb F_p$. Taking $y,y' \in
\cohp 1 \bq$ such that $\delta_p(y)=\alpha'$, and $\delta_p(y')=\beta'$, we see
that $$\phi(\gamma')=\delta_p\pone\pi_*(\gamma')=
\delta_p\pone(\beta_p(y)y'-\beta_p(y')y)=
\delta_p(\beta_p(y)^py'-\beta_p(y')^py)=\alpha'^p\beta'-\beta'^p\alpha'.$$
It follows that
$$\alpha^p\beta-\beta^p\alpha=\inf(\alpha'^p\beta'-\beta'^p\alpha')=
\inf(\phi(\gamma'))=\phi\inf(\gamma')=0.$$ \qed\par
We note that the result on $d_{2p-1}$ could also
be considered as a case of an integral version
of Kudo's transgression theorem. We are now ready to state our theorem on
$\co * \tilp$.
\proclaim Theorem 2. Let $p$ be an odd prime, and let $\tilp$ be the
group defined above. Then $H^*(\tilp;\Bbb Z)$ is generated by elements
$\alpha$,\allowbreak$\beta$,\allowbreak$\chi_1,\ldots,\chi_{p-1}$,
\allowbreak$\zeta$, with
$$\deg(\alpha)=\deg(\beta)=2 \quad \deg(\chi_i)=2i \quad \deg(\zeta)=2p,$$
subject to the following relations:
$$p\alpha=p\beta=0$$
$$\alpha^p\beta=\beta^p\alpha$$
\hfil\hbox{$\alpha\chi_i=\cases{0 \cr -\alpha^p}\qquad
\beta\chi_i=\cases{0 &for $i<p-1$ \cr -\beta^p &for $i=p-1$}$}\hfil
$$\chi_i\chi_j=\cases{p\chi_{i+j} &$i+j<p$\cr p^2\zeta &$i+j=p$ \cr
p\zeta\chi_{i+j-p} &$p<i+j<2p-2$ \cr
p\zeta\chi_{p-2}+\alpha^{2p-2}+\beta^{2p-2}-\alpha^{p-1}\beta^{p-1}
&$i=j=p-1$}$$ \pra
Chern classes of representations of $\tilp $ generate the whole ring.
An automorphism of $\tilp$ sends $\chi_i$ to $\chi_i$
(resp.\ $(-1)^i\chi_i$) and $\zeta$ to $\zeta$ (resp.\ $-\zeta$) if it
fixes (resp. reverses) \sone. The effect of an automorphism on $\alpha$,
$\beta$ may be determined from their definition. Considered as elements of
$\hom(\tilp,\sone)$, $\alpha$ has kernel $\langle \sone, B\rangle$ and sends
$A$ to $e^{2\pi i/p}$, and $\beta$ has kernel $\langle \sone, A\rangle$ and
sends $B$ to $e^{2\pi i/p}$.
If we let $H$ be the subgroup generated by $B$ and elements of \sone\
we may define $$\chi_i=\cases{\cor_H^{\tilp}(\tau'^i) &for $i<p-1$ \cr
\cor_H^{\tilp}(\tau'^{p-1})-\alpha^{p-1} &for $i=p-1$}$$
where $\tau'$ is any element of $\cohz 2 H $ restricting to \sone\ as the
generator $\tau$.
Similarly, $\zeta =c_p(\rho)$, where $\rho$ is an irreducible representation
of $\tilp$ restricting to \sone\ as $p$ copies of the representation
$\xi$ with $c_1(\xi)=\tau $. \par
\proof First we note that in the $E_\infty$ page of the above spectral sequence
all the group extensions that we need to examine are extensions of finite
groups by the infinite cyclic group, so are split. The elements $\alpha$ and
$\beta$ defined in the statement above clearly yield generators for
$E_\infty^{2,0}$, and the relations between them are exactly the relations that
hold between the corresponding elements in the spectral sequence.
Let $\beta' $ in $\co 2 {{\rm B}H}$ be the restriction to $H$ of
$\beta$, and take any choice of $\tau'$ as in the statement. We may show by
considering $\beta'$ and $\tau'$ as homomorphisms from $H$ to \sone\ that
conjugation by $A^i$ induces the map on $\co 2 {{\rm B}H}$ that fixes $\beta'$
and sends $\tau'$ to $\tau'-i\beta'$. Now applying the formula for
$\res^G_K\cor^G_H$ (see for example [\brown]) it follows that $\chi_i$
restricts to \sone\ as $p\tau^i$, so yields a generator for $E_\infty^{0,2i}$.
\par
Any irreducible representation of $\tilp$ has degree 1 or $p$, because $\tilp$
has an abelian subgroup of index $p$.
Let $\rho$ be the representation of $\tilp$ induced from a 1-dimensional
representation of $H$ with first Chern class $\tau'$. $\rho$ restricts to
\sone\ as $p$ copies of the representation with first Chern class $\tau$, so
its total Chern class restricts to \sone\ as $(1+\tau)^p$, and so $c_p(\rho)$
yields a generator for $E_\infty^{0,2p}$, and
$c_i(\rho)={1/p}{p\choose i}\chi_i+P_i(\alpha,\beta)$ for some polynomial
$P_i$. We shall show later that $P_i =0$.
\par
The restriction to $H$ of $\alpha$ is trivial, so by Fr\"obenius reciprocity
$$\alpha\cor^\tilp_H(\tau'^i)=\cor^\tilp_H(\res^\tilp_H(\alpha)\tau'^i)=0,$$
and the expressions given for $\alpha\chi_i$ follow. By calculating
$\alpha(\beta\chi_i)=\beta(\alpha\chi_i)$, we may deduce that $\beta\chi_i=0$
for $i < p-1$ and $\beta\chi_{p-1}=\lambda(\alpha^{p-1}\beta-\beta^p)
-\alpha^{p-1}\beta$ for some scalar $\lambda$. To show that $\lambda =1$ we
use the restriction map to $H$, and the formula for corestriction followed by
restriction.
$$\eqalign{\res_H^\tilp(\beta \chi_{p-1})&=\beta^\prime \sum_{i=0}^{p-1}
(\tau^\prime+i\beta^\prime)^{p-1}\cr
&= \beta^\prime \sum_{j=0}^{p-1} \tau^{\prime p-1-j} \beta^{\prime j}
\sum _{i=0}^{p-1} i^j}$$
Newton's formula tells us that
$$\sum_{i=1}^{p-1} i^j \equiv \cases{0 \enspace (p) &for $j \not\equiv 0
\enspace (p-1)$ \cr 1\enspace (p) &for $j \equiv 0 \enspace (p-1)$ }$$
so $\res^\tilp_H(\beta\chi_{p-1})=-\beta'^p$, and the required relation
follows. \par
We now know $\res^\tilp_\sone(\chi_i\chi_j)$, $\alpha\chi_i\chi_j$, and
$\beta\chi_i\chi_j$, which together imply the relations given for
$\chi_i\chi_j$. To complete the proof of the theorem we must determine the
effect of automorphisms of $\tilp$ on the $\chi_i$. We know that an
automorphism sends $c_i(\rho)$ to itself or $(-1)^i$ times itself depending
whether or not it reverses the sense of \sone, so it will suffice to show that
$\chi_i={1/p}{p\choose i}c_i(\rho)$. The character of $\rho$ is zero except on
$\sone$, so if $\theta $ is a 1-dimensional representation of $\tilp$
restricting trivially to \sone, then $\rho\otimes\theta$ is isomorphic to
$\rho$. If we apply the formula expressing $c.(\rho\otimes\theta)$ in terms of
$c.(\rho)$ and $c.(\theta)$ (see [\atiyah]) we obtain
$$c_i(\rho)=c_i(\rho \otimes \theta)=
\sum_{j=0}^i {p-i+j \choose j}c_1(\theta)^j c_{i-j}(\rho).$$
and hence inductively
$$c_i(\rho)c_1(\theta)=\cases {0 &for $i<p-1$ \cr -c_1(\theta)^p &for $i=p$}.$$
Since $\alpha$ and $\beta$ are possible values for $c_1(\theta)$ the required
result follows. We may show inductively that $\chi_i$ is in the subring
generated by Chern classes because $\chi_1$ is, and $\chi_1\chi_{i-1}$,
${1/p}{p\choose i}\chi_i$ are coprime multiples of $\chi_i$.
\qed

We are now ready to state our theorem on the integral cohomology of ${\rm
B}P(n)$.
\proclaim Theorem 3. Let p be an odd prime and let $P(n)$ be as
defined above. Then $H^*({\rm B}P(n);\Bbb Z)$ is generated by elements
$\alpha,\allowbreak\beta,\allowbreak\mu,\allowbreak\nu,\allowbreak
\chi_1,\ldots,\chi_{p-1}, \allowbreak\zeta$, with
$$\deg(\alpha)=\deg(\beta)=2\quad \deg(\mu)=\deg(\nu)=3\quad \deg(\chi_i)=2i
\quad\deg(\zeta)=2p$$
subject to the following relations:
$$p\alpha=p\beta=0\quad p\mu=p\nu=0\quad p^{n-3}\chi_1=0\quad p^{n-2}\chi_i=0
\quad p^{n-1}\zeta=0$$
$$\alpha\mu=\beta\nu$$
$$\alpha^p\beta=\beta^p\alpha\quad \alpha^p\mu=\beta^p\nu$$
$$\hbox{$\alpha\chi_i=\cases{0 & \cr -\alpha^p & }$\qquad
$\beta\chi_i=\cases{0 &for $i<p-1$ \cr -\beta^p &for $i=p-1$}$}$$
$$\hbox{$\mu\chi_i=\cases{0 & \cr -\beta^{p-1}\mu & }$\qquad
$\nu\chi_i=\cases{0 &for $i<p-1$ \cr -\alpha^{p-1}\nu &for $i=p-1$}$}$$
$$\chi_i\chi_j=\cases{p\chi_{i+j} &$i+j<p$ \cr
p^2\zeta &$i+j=p$ \cr p\zeta\chi_{i+j-p} &$p<i+j<2p-2$ \cr
p\zeta\chi_{p-2}+\alpha^{2p-2}+\beta^{2p-2}-\alpha^{p-1}\beta^{p-1}
&$i=j=p-1$}$$
$$\mu\nu=\cases{0 &for $n>3$ \cr \lambda\chi_3 &for $n=3,\ p>3, \lambda \in
\Bbb Z_p^\times$ \cr
3\lambda\zeta &for $n=3,\ p=3,\  \lambda=\pm 1$}$$
\pra
Chern classes of representations of $P(n)$ generate $H^{\rm even}({\rm B}P(n);
\Bbb Z)$.
Under an automorphism of $P(n)$ which restricts to the centre as
$C \mapsto C^j$, $\chi_i$ is mapped to $j^i\chi_i$, and $\zeta$ is mapped to
$j^p\zeta$. The effect of automorphisms on $\alpha$ and $\beta$ is
determined by the natural isomorphism $H^2({\rm B}G;\Bbb Z) \cong
\hom(G,\Bbb R/\Bbb Z)$, under which \pra
\hfil\hbox{$\eqalign{\alpha :A &\mapsto 1/p \cr
B &\mapsto 0 \cr C &\mapsto 0 }\quad \eqalign{\beta : A &\mapsto 0 \cr
B &\mapsto 1/p \cr C &\mapsto 0} \quad \eqalign{\chi_1 : A &\mapsto 0 \cr
B &\mapsto 0\cr C &\mapsto 1/{p^{n-3}}.}$}\hfil
\pra
An automorphism of $P(n)$ which sends $\alpha$ to $n_1\alpha+n_2\beta$,
$\beta$ to $n_3\alpha+n_4\beta$ and restricts to the centre as $C\mapsto C^j$
sends $\mu$ to $j(n_4\mu+n_3\nu)$ and $\nu$ to $j(n_2\mu+n_1\nu)$.
If $\gamma'$ in $H^2({\rm B}\langle B,C \rangle;\Bbb Z)$ is such that it
maps to the
following element of $\hom(\langle B,C \rangle, \Bbb R/\Bbb Z)$
$$\eqalign{\gamma ' : B &\mapsto 0 \cr C &\mapsto 1/{p^{n-2}}},$$
then $\chi_i$ is defined as follows:
$$\chi_i=\cases {\cor _{\langle B,C \rangle}^{P(n)}(\gamma'^i) &for $i<p-1$
\cr \cor_{\langle B,C \rangle}^{P(n)}(\gamma'^{p-1})-\alpha^{p-1}
&for $i=p-1$.}$$
\pra
These are, up to scalar multiples, equal to $c_i(\rho)$, where $\rho$ is a
$p$-dimensional irreducible representation of $P(n)$, whose restriction
to $\langle C \rangle$ is a sum of $p$ copies of the representation $\theta$,
with $c_1(\theta)=\res_{\langle C\rangle}^{\langle B,C \rangle}(\gamma')$.
In fact,
$ c_i(\rho)=\textstyle{{1 / p} {p \choose i}}\chi_i.$
Also, we may define $\zeta= c_p(\rho)$.
\par
\proof We examine the spectral sequence for ${\rm B}P(n)$ as an \sone-bundle
over ${\rm B}\tilp$. $E_2^{*,0}$ is isomorphic to $\cohz * \bpn$ and
$E_2^{*,*}$ is freely generated by $E_2^{*,0}$ and an element $\xi$ of infinite
order in $E_2^{0,1}$. We know that $\co 2 {{{\rm B}P(n)}}\cong\hom(P(n),\sone)
\cong C_{p^{n-3}}\oplus C_p\oplus C_p$, so $d_2(\xi)$ must be
$\pm p^{n-3}\chi_1$. If we wanted to calculate the cohomology of the metacyclic
groups $M(n)$ described above, the differential
in this spectral sequence would send $\xi$ to
$\pm p^{n-3}\chi_1+\gamma$ for some non-zero $\gamma$ in $\langle \alpha,\beta
\rangle$. It is now easy to see that $E_\infty$ is generated by the elements
$\alpha,\allowbreak\beta,\allowbreak
\mu=\beta\xi,\allowbreak\nu=\alpha\xi,\allowbreak
\chi_1,\ldots,\chi_{p-1}$ and $\zeta$ subject to the relations they satisfy as
elements of $E_2^{*,*}$ together with $p^{n-3}\chi_1=0$, $p^{n-2}\chi_i=0$, and
$p^{n-1}\zeta=0$. For each~$m$, the filtration of $\co m {{\rm B}P(n)}$ given by
the $E_\infty$ page is trivial, so we may use the same symbols to denote
elements of $\co m {{\rm B}P(n)}$, and the relations that hold in $E_\infty$
determine all the relations that hold in $\co m {{\rm B}P(n)}$ except for the
product of the two odd dimensional generators.
\par
We know that $p\mu\nu=0$, and
the relation $\alpha\mu=\beta\nu$ implies that $\alpha\mu\nu=\beta\mu\nu=0$,
and so $\mu\nu$ must be a multiple of $p^{n-3}\chi_3$ for $p\geq 5$ (resp.\
$3\zeta$ for $p=3$). Note that these elements restrict to zero on all proper
subgroups of $P(n)$.
In the case of $P(3)$, Lewis [\lew] shows that $\mu\nu$ is not zero by
considering the spectral sequence for $P(n)$ considered as an extension of a
maximal subgroup by $C_p$. A similar method will work in general, but
we offer an alternative proof that involves expressing $\mu$ and
$\nu$ as Bocksteins of elements of $\cohp 2 {{\rm B}P(n)}$. This proof is
contained in lemma~4 and corollary~5.
\par
The effect of automorphisms on $\chi_i$ and $\zeta$ is easily seen to be as
claimed from their alternative definitions as Chern classes. To determine the
effect of automorphisms on $\mu$ and $\nu$, we note that an automorphism of
$P(n)$ restricting to the centre as $C\mapsto C^j$ extends to an endomorphism
of $\tilp$ which wraps the central circle $j$ times around itself, so induces a
map of the above spectral sequence to itself sending $\xi$ to $j\xi$. This
completes the proof of theorem~3 modulo lemma~4 and its corollary.\qed\par

We now examine the spectral sequence with $\Bbb F_p$ coefficients for the
central extension $C_{p^{n-2}}\rightarrowtail P(n)\twoheadrightarrow C_p\oplus
C_p$. Take generators so that $\cohp * {{\rm B}C_p\oplus C_p}\cong
\Bbb F_p[x,x']\otimes \Lambda[y,y']$, where $\beta_p(y)=x$, $\beta_p(y')=x'$,
and $\cohp * {{\rm B}C_{p^{n-2}}}\cong \Bbb F_p[t]\otimes\Lambda[u]$, where
$\beta_p(u)=t$ for $n=3$ (resp. $\beta_p(u)=0$ for $n\geq 4$). Then the $E_2$
page is isomorphic to $\Bbb F_p[x,x',t]\otimes \Lambda[y,y',u]$, and the first
two differentials are as described in the following lemma.
\proclaim Lemma 4. With notation as above, identify the elements $x,x',y,y'$
in the spectral sequence
with their images in $\cohp * {{\rm B}P(n)}$ under the inflation map.
\item{1)} Let $n\geq 4$. Then $d_2$ is trivial, and $d_3(t)$ is a non-zero
multiple of $xy'-x'y$. The set $\{x,x',yy',u'y,u'y'\}$ is a basis for $\co 2
\bpn$, where $u'$ is any element of $\co 1 \bpn$ restricting to $C_{p^{n-2}}$
as $u$. \pra
\item{2)} Let $n=3$. Then $d_2(u)$ is a non-zero multiple of $yy'$, $d_2(t)=0$,
and $E_3$ is generated by $y,y',x,x',[uy],[uy']$ and $t$ subject to the
relation $yy'=0$ and those implied by the relations in $E_2$. In particular
$[uy]y'=-[uy']y$ but this element is non-zero. As in the case $n\geq 4$,
$d_3(t)$ is a non-zero multiple of $xy'-x'y$. Let $Y,Y'$ be elements of $\co 2
\bpn$ such that $\{x,x',Y,Y'\}$ is a basis for $\co 2 \bpt$, and let
$X=\beta_p(Y)$, $X'=\beta_p(Y')$. Then $\{yY',xy,xy',x'y',X,X'\}$ is a basis
for $\co 3 \bpt$ and $\{xX,xX',x'X',x^2y,x^2y',xx'y',x'^2y',YX'\}$ is a basis
for $\co 5 \bpt$.
\par
\proof 1). In this case $H^1$ has order $p^3$, so $u$ must survive. The element
$xy'-x'y$ is the image under $\pi_*$ of a generator for $\cohz 3 {{\rm
B}(C_p\oplus C_p)}$, so must be killed by some differential. We have already
shown that it cannot be killed by $d_2$, so the only possibility is that $t$
survives until $E_3$ and kills it. The rest of the statement follows easily.
\par
2). In this case $H^1$ has order $p^2$, so $d_2(u)$ must be non-zero. It is
true in general that if $G$ is a central extension of $C_p$ by $Q$, then in the
corresponding spectral sequence with $\Bbb F_p$ coefficients
$d_2:E_2^{0,1}\rightarrow E_2^{2,0}$ must kill the extension class. This
follows by naturality, since one may regard the extension class as defining a
homotopy class of maps from $\bq$ to $K(C_p,2)$ such that $\bg$ is the
 ${\rm B}C_p$-bundle induced by the path-loop fibration over $K(C_p,2)$. Since
all subgroups of $P(3)$ of order $p^2$ are copies of $C_p\oplus C_p$, the
extension class of $P(3)$ must restrict to zero on all cyclic subgroups, so
must be a multiple of $yy'$. The transgression commutes with the Bockstein so
$d_2(t)=0$ and $d_3(t)=\beta_pd_2(u)$.
\par
Given the values of these differentials it is routine to compute the $E_4$ page
of the spectral sequence. If we write $E_r^n=\bigoplus_{i+j=n}E_r^{i,j}$, then
$\{[uy],[uy'],x,x'\}$ forms a basis for $E_4^2=E_\infty^2$, and
$\{[ty],[ty'],[uy]y',xy,xy',x'y'\}$ forms a basis for $E_4^3=E_\infty^3$. The
spectral sequence operation $_F\beta$ introduced by Araki [\araki] and Vasquez
[\vasquez] maps $[uy]$ to $[ty]$ and $[uy']$ to $[ty']$, so if $Y$ and $Y'$ are
chosen to yield the generators for $E_4^{1,1}$ their Bocksteins yield
generators for $E_4^{1,2}$. A basis for $E_4^5$ is given by the eight elements
of the statement, which we know to be universal cycles, and the elements
$[t^2y]$, $[t^2y']$. $E_4^4$ consists of universal cycles, and the universal
coefficient theorem tells us that $H^5$ has order $p^8$, so
$[t^2y]$ and $[t^2y']$ cannot be universal cycles. \qed \par
\proclaim Corollary 5. In $\cohz * \bpn $ the product $\mu\nu$ is non-zero if
and only if $n=3$. \par
\proof In the notation of lemma~4 it suffices to determine
$\delta_p(u'y)\delta_p(u'y')$ in the case $n\geq 4$, and
$\delta_p(Y)\delta_p(Y')$ in the case $n=3$. In the case when $n=3$,
$$\delta_p(Y)\delta_p(Y')=\delta_p(Y\beta_p(Y'))=\delta_p(YX').$$
The kernel of $\delta_p : \cohp 5 \bpt \rightarrow \cohz 6 \bpt$ is equal to
$\pi_*(\cohz 5 \bpt)$, which is generated by $xX$, $xX'$ and $x'X'$, so by
lemma~4 $\delta_p(YX')$ is non-zero.\par
In the case when $n=4$, $\cohz i \bpf $ has exponent $p$ for $i=2,3$, so
$\pi_*$ is injective from these groups, and $\ker\beta_p:\co 2 \bpf \rightarrow
\co 3 \bpf$ is equal to $\beta_p(\co 1 \bpf)$. $\beta_p(yy')=xy'-x'y=0$, so we
may choose the element $u'$ in lemma~4 so that $\beta_p(u')=\lambda yy'$ for
some non-zero $\lambda$. Then we have
$$\delta_p(u'y)\delta_p(u'y')=\delta_p(u'y\beta_p(u'y'))=
\delta_p(u'y(\lambda yy'y'-u'x'))=0.$$
The case when $n\geq 5$ is similar but simpler, since $u'$ may be chosen so
that $\delta_p(u')=p^{n-4}\chi_1$, which implies that $\beta_p(u')=0$. \qed
\par
\noindent{\bf Remarks.} Theorem~3 contains independent proofs of Thomas'
result that the even degree subring of
$\cohz * \bpn$ is generated by Chern classes [\thompap],
and Lewis' calculation of $\cohz * \bpt $. Our notation differs slightly from
that of Lewis. We have renumbered the generators $\chi_i$ (note that
$\chi_1$ vanishes
for $n=3$). Also our $\chi_{p-1}$ and Lewis' $\chi_{p-2}$ are related
by the formula $$\chi_{p-2}^{\rm Lewis}=\chi_{p-1}+\alpha^{p-1}+\beta^{p-1}.$$
Our result disagrees with that of AlZubaidy [\alzub].
\par
The author wishes to thank his supervisor C.~B.~Thomas for encouragement with
this work, P.~H.~Kropholler for suggesting the Lie group construction, and
J.~Huebschmann for helpful comments concerning an earlier version of this work.
\par

\beginsection{{\rm REFERENCES}}\par
\paper \alzub/K.~A/L{\rm Z}UBAIDY/Rank 2 $p$-groups, $p>3$, and Chern classes/
Pacific J. Math./103/1982/259--267/
\paper \araki/S.~A/RAKI/Steenrod reduced powers in the spectral sequences
associated with a fibering I, II/Mem. Fac. Sci. Kyusyu Univ. Series (A) Math./
11/1957/15--64 81--97/

\paper \atiyah/M.~F.~A/TIYAH/Characters and the cohomology of finite groups/
Publ. Math. IHES/9/1961/23--64/
\book \brown/K.~S.~B/ROWN/Cohomology of groups/Springer Verlag/1982/

\prepaper \hueb/J.~H/UEBSCHMANN/Perturbation theory and free resolutions for
nilpotent groups of class 2/{\it J. of Algebra} (to appear)./
\prepaper \huec/J.~H/UEBSCHMANN/Cohomology of nilpotent groups of class 2/
{\it J. of Algebra} (to appear)./
\prepaper \leary/I.~J.~L/EARY/The mod-$p$ cohomology rings of some $p$-groups/
In preparation./
\paper \lew/G.~L/EWIS/The integral cohomology rings of groups of order $p^3$/
Trans. Amer. Math. Soc./132/1968/501--529/
\prepaper \moselle/B.~M/OSELLE/Calculations in the cohomology of finite groups/
Unpublished essay, (1988)./
\paper \thompap/C.~B.~T/HOMAS/Riemann-Roch formulae for group representations/
Mathematika/20/1973/253--262/
\book \Thomas/C.~B.~T/HOMAS/Characteristic classes and the cohomology of finite
groups/Cambridge University Press/1986/
\paper \vasquez/R.~V/ASQUEZ/Nota sobre los cuadrados de Steenrod en la sucesion
espectral de un espacio fibrado/Bol. Soc. Mat. Mexicana/2/1957/1--8/

\end